\magnification=1200
\newif\ifdraft\draftfalse
\ifdraft\input colordvi\else\input blackdvi\fi%
\overfullrule=0pt

\global \edef \lmlbvjac {\Red {1.2}}
\global \edef \lemmonomrep {\Red {1.3}}
\global \edef \thmintclmonom {\Red {2.3}}
\global \edef \propmonomialval {\Red {3.2}}
\global \edef \corcorr {\Red {3.3}}
\global \edef \lmadjmonom {\Red {3.4}}
\global \edef \sectgenmonRVa {\Red {4}}
\global \edef \thmadjmonomval {\Red {4.1}}
\global \edef \corsubaddmon {\Red {4.2}}
\global \edef \sectchar {\Red {5}}
\global \edef \propreesvalcon {\Red {5.1}}
\global \edef \rmkreesval {\Red {5.1}}
\global \edef \thmadjone {\Red {5.2}}
\global \edef \rmkpseudorat {\Red {5.2}}
\global \edef \rmkreesvaluni {\Red {5.2}}
\global \edef \cormonompsrat {\Red {5.3}}
\global \edef \coradjtwo {\Red {5.4}}
\global \edef \rmksubaddtwo {\Red {5.4}}
\expandafter \gdef \csname CHRR\endcsname {1}
\expandafter \gdef \csname DEL\endcsname {2}
\expandafter \gdef \csname Howald\endcsname {3}
\expandafter \gdef \csname Hu87\endcsname {4}
\expandafter \gdef \csname HunSwa\endcsname {5}
\expandafter \gdef \csname Kap\endcsname {6}
\expandafter \gdef \csname KiSt\endcsname {7}
\expandafter \gdef \csname LipmanRS\endcsname {8}
\expandafter \gdef \csname Lipmanadj\endcsname {9}
\expandafter \gdef \csname Lipproxy\endcsname {10}
\expandafter \gdef \csname LipSat\endcsname {11}
\expandafter \gdef \csname LipTei\endcsname {12}
\expandafter \gdef \csname MS63\endcsname {13}
\expandafter \gdef \csname Ra74\endcsname {14}
\expandafter \gdef \csname Rees\endcsname {15}
\expandafter \gdef \csname SwHu\endcsname {16}
\expandafter \gdef \csname TakWat04a\endcsname {17}

\headline={}
\headline={\ifnum\pageno>1\ifodd\pageno{\tenrm Adjoints of 
ideals\hfill\the\pageno}\else{\the\pageno\tenrm\hfill R. H\"ubl and I. Swanson}%
	\fi\fi}
\footline={}
\ifdraft\footline={\hfill\today}\fi

\font\footfont=cmr8
\font\footitfont=cmti8

\font\largerm=cmr12 at 15pt
\font\medi=cmr12 at 12 pt
\font\bb=msbm10
\font\Cal=eusm10 at 10pt 
\font\sCal=eusm7
\font\goth=eufm10
\font\sgoth=eufm7
\font\ssgoth=eufm5
\font\ftdefit=cmmib10 
\font\ftdefits=cmmib7 

\def\df{\bf}
\def\bfi#1{\hbox{\ftdefit#1}} 
\def\bfis#1{\hbox{\ftdefits#1}} 

\def\adj{\mathop{\rm adj}\nolimits}

\def\bbQ{{\hbox{\bb Q}}}

\def\m{{\hbox{\goth m}}}
\def\p{{\hbox{\goth p}}}
\def\P{{\hbox{\goth P}}}

\def\q{{\hbox{\goth q}}}
\def\V{{\hbox{\goth V}}}
\def\n{{\hbox{\goth n}}}

\def\sp{{\hbox{\sgoth p}}}
\def\sq{{\hbox{\sgoth q}}}
\def\ssp{{\hbox{\ssgoth p}}}
\def\sP{{\hbox{\sgoth P}}}
\def\sV{{\hbox{\sgoth V}}}

\def\inc{\subseteq}
\def\ord{\mathop{\rm ord}\nolimits}
\def\ic#1{\,\overline{\!#1}\,\!}
\def\bbN{\hbox{\bb N}}
\def\bbZ{\hbox{\bb Z}}

\def\bbQ{\hbox{\bb Q}}
\def\Reesval{\hbox{\Cal RV}}
\def\sReesval{\hbox{\sCal RV}}
\def\val{\hbox{\Cal V}}
\def\sval{\hbox{\sCal V}}
\def\NP{\hbox{NP}}
\def\NPint{\hbox{NP}^{\,\circ}}
\def\th{{$^{\hbox{\sevenrm th}}$}}

\def\qedbox{\hbox{\vbox{\hrule\hbox{\vrule\kern3pt\vbox{\kern6pt}\kern3pt\vrule}\hrule}}}

\def\qed{\hfill\qedbox\vskip1.1ex} 

\let\text\hbox
\def\overset#1\to#2{\mathop{\buildrel #1 \over #2}}
\def\overarrow#1{\mathop{\buildrel #1 \over \rightarrow}}
\def\frac#1#2{{#1 \over #2}}
\font\tenmsb=msbm10
\font\sevenmsb=msbm7
\font\fivemsb=msbm5
\newfam\msbfam
\textfont\msbfam=\tenmsb
\scriptfont\msbfam=\sevenmsb
\scriptscriptfont\msbfam=\fivemsb

\def\hexnumber#1{\ifcase#1 0\or1\or2\or3\or4\or5\or6\or7\or8\or9\or
	A\or B\or C\or D\or E\or F\fi}

\mathchardef\subsetneq="2\hexnumber\msbfam28

\def\optioncite[#1]#2{[\csname #2\endcsname, #1]}
\def\nooptioncite#1{[\csname #1\endcsname]}
\def\cite{\futurelet\testchar\maybeoptioncite}
\def\maybeoptioncite{\ifx[\testchar \let\next\optioncite
	\else \let\next\nooptioncite \fi \next}

\newwrite\isauxout
\openin1\jobname.aux
\ifeof1\message{No file \jobname.aux}
       \else\closein1\relax\input\jobname.aux
       \fi
\immediate\openout\isauxout=\jobname.aux
\immediate\write\isauxout{\relax}
\def\draftcmt#1{\llap{\smash{\raise2ex\hbox{{#1}}}}}
\font\foottt=cmtt8
\def\label#1{{\global\edef#1{\the\sectno.\the\thmno}}}
\def\pagelabel#1{{\global\edef#1{\the\pageno}}}
\def\label#1{\ifdraft\draftcmt{\foottt\Red{\string#1}\ \ }\unskip\fi%
	\immediate\write\isauxout{\noexpand\global\noexpand\edef\noexpand#1%
		{\noexpand\Red{\the\sectno.\the\thmno}}}%
	{\global\edef#1{\noexpand\Red{\the\sectno.\the\thmno}}}%
	\unskip\ignorespaces}

\def\today{\ifcase\month\or January\or February\or March\or
  April\or May\or June\or July\or August\or September\or
  October\or November\or December\fi
  \space\number\day, \number\year}

\def\widow#1{\vskip 0pt plus#1\vsize\goodbreak\vskip 0pt plus-#1\vsize}

\newcount\sectno \sectno=0
\newcount\thmno \thmno=0

\newdimen\sectskip \newdimen\secttitleskip
\sectskip=4ex \secttitleskip=3ex
\def\sectionbreak{\penalty-9999}%

\def\labelsection[#1]#2{%
  \noindent{\bgroup\def\noexpand{\relax}%
          \bf\Maroon{\the\sectno. #2}\egroup\nobreak%
          \ifdraft\smash{\rlap{\foottt\Red{\string#1\ \ }}}\fi}%
  \immediate\write\isauxout{\noexpand\global\noexpand\edef\noexpand#1%
                {\noexpand\Red{\the\sectno}}}%
  \global\edef#1{{\the\sectno}}%
  \vskip\secttitleskip%
  \hrule height0pt
  \unskip%
}

\def\nolabelsection#1{%
  \noindent{\bgroup\def\noexpand{\relax}%
          \bf\Maroon{\the\sectno. #1}\egroup}\nobreak%
	\vskip\secttitleskip%
	\hrule height0pt
	\unskip%
  }

\def\section{\widow{.05}
	\removelastskip\global\advance\sectno by1%
	\global\thmno=0\sectionbreak\vskip\sectskip%
	\futurelet\testchar\maybelabelsection}
\def\maybelabelsection{\ifx[\testchar\let\next\labelsection%
	\else\let\next\nolabelsection\fi\next%
}
\def \subsection#1{\vskip 0.5truecm\relax
	\noindent{\bf #1} \vskip 0.4truecm\relax}
\def \thmline#1{\vskip .2cm\relax
	\global\advance\thmno by 1
	\noindent{\bf #1\ \the\sectno.\the\thmno:}\ \ %
	\bgroup \advance\baselineskip by -1pt \it
	\abovedisplayskip =4pt
	\belowdisplayskip =3pt
	\parskip=0pt
	}
\def \wthmline#1{\vskip .2cm\relax
	\noindent{\bf #1:}\ \ %
	\bgroup \advance\baselineskip by -1pt \it
	\abovedisplayskip =4pt
	\belowdisplayskip =3pt
	\parskip=0pt
	}
\def \dthmline#1{\vskip 6pt\relax
	\noindent{\bf #1:}\ \ }
\def \defin{\thmline{Definition}}
\def \thm{\thmline{Theorem}}

\def \endb{\egroup \vskip 0.2cm\relax}
\def \cor{\thmline{Corollary}}
\def \lemma{\thmline{Lemma}}
\def \remark{\dthmline{Remark}}

\def \example{\dthmline{Example}}
\def \prop{\thmline{Proposition}}
\def \proof{\vskip0.0cm\relax\noindent {\sl Proof:\ \ }}
\def \proofof#1{\medskip\noindent {\sl Proof of #1:\ }}

\ %
\vskip 3ex
\centerline{\largerm  Adjoints of ideals}
\vskip 2ex
\centerline{\medi Reinhold H\"ubl
and Irena Swanson\footnote{$^*$}{\footitfont Partially supported by the National
Science Foundation}}
\unskip\footnote{ }{{\footitfont 2000 Mathematics Subject Classification}
\footfont 13A18, 13A30, 13B22, 13H05.}
\unskip\footnote{ }{{\footitfont Key words and phrases.}
\footfont Adjoint ideal, multiplier ideal, regular ring, Rees valuation.}

\bigskip

{\bgroup
\baselineskip=10pt
\narrower\narrower
\tenrm
\noindent
{\tenbf Abstract.}
We characterize ideals whose adjoints are determined
by their Rees valuations.
We generalize the notion of a regular system of parameters,
and prove that for ideals generated by monomials in such elements,
the integral closure and adjoints are generated by monomials.
We prove that the adjoints of such ideals
and of all ideals in two-dimensional regular local rings
are determined by their Rees valuations.
We prove special cases of subadditivity of adjoints.

\egroup}

\bigskip

Adjoint ideals and multiplier ideals have recently emerged
as a fundamental tool in commutative algebra and algebraic geometry.
In characteristic 0 they may be defined using resolution of singularities.
In all characteristics, even mixed,
Lipman gave the following definition:

\defin
Let $R$ be a regular domain,
$I$ an ideal in $R$.
The {\bf adjoint}  $\adj I$ of $I$ is defined as follows:
$$
\adj I = \bigcap_v\, \{r \in R ~|~ v(r) \ge v(I) - v(J_{R_v/R}) \},
$$
where the intersection varies over all valuations $v$
on the field of fractions $K$ of $R$ that are non-negative on $R$
and for which the corresponding valuation ring $R_v$
is a localization of a finitely generated $R$-algebra.
The symbol $J_{R_v/R}$ denotes the Jacobian ideal of $R_v$ over~$R$.
\endb

By the assumption on $v$,
each valuation in the definition of $\adj I$ is Noetherian.

Many valuations $v$ have the same valuation ring $R_v$;
any two such valuations are positive real multiples of each other,
and are called equivalent.
In the definition of $\adj I$ above,
one need only use one $v$ from each equivalence class.
In the sequel,
we will always choose {\bf normalized} valuations,
that is,
the integer-valued valuation $v$ such that for all $r \in R$,
$v(r)$ equals that non-negative integer $n$
which satisfies that
$r R_v$ equals the $n$th power of the maximal ideal of $R_v$.

Lipman proved that for any ideal $I$ in $R$ and any $x \in R$,
$\adj(xI) = x \adj(I)$.
In particular,
$\adj(xR) = (x)$.

A crucial and very powerful property is the subadditivity of adjoints:
$\adj(IJ) \inc \adj(I) \adj(J)$.
This was proved in characteristic zero
by Demailly, Ein and Lazarsfeld~\cite{DEL},
and is unknown in general.
We prove it for generalized monomial ideals in Section~\sectgenmonRVa,
and for ideals in two-dimensional regular domains in Section~{\sectchar}.
The case of subadditivity of adjoints for ordinary monomial ideals can be deduced from
Howald's work~\cite{Howald} using toric resolutions,
and the two-dimensional case has been proved by
Takagi and Watanabe~\cite{TakWat04a} using multiplier ideals.
The case for generalized monomial ideals proved here is new.

An aspect of proving subadditivity and computability of adjoints
is whether there are
only finitely many valuations $v_1, \ldots, v_m$ such that for all $n$,
$$
\adj (I^n)
= \bigcap_{i=1}^m \{r \in R ~|~ v_i(r) \ge v_i(I^n) - v_i(J_{R_{v_i}/R}) \}.
$$
We prove in Sections~{\sectgenmonRVa}
that Rees valuations suffice for the generalized monomial ideals.
We also give an example (the first example in Section~\sectchar)
showing that Rees valuations do not suffice in general.
In Section~{\sectchar} we give a general criterion
for when the adjoint of an ideal is determined by its Rees valuations.
A corollary is that Rees valuations suffice for
ideals in two-dimensional regular domains.
The first three sections develop the background on
generalized monomial ideals.

\section{Generalized regular system of parameters}

\defin
Let $R$ be a regular domain.
Elements $x_1, \ldots, x_d$ in $R$
are called a {\bf generalized regular system of parameters} if
$x_1, \ldots, x_d$ is a permutable regular sequence in $R$
such that for every $i_1, \ldots, i_s \in \{1, \ldots, d\}$,
$R/(x_{i_1}, \ldots, x_{i_s})$ is a regular domain.
\endb

\remark
Any part of a generalized regular system of parameters is
again a generalized regular system of parameters.

\medskip

For example, when $R$ is regular local,
an arbitrary regular system of parameters
(or a part thereof) is a generalized regular system of parameters;
or if $R$ is a polynomial ring over a field,
the variables are a generalized regular system of parameters.

Let $\p$ be any prime ideal
containing the generalized regular system of parameters $x_1, \ldots, x_d$.
As $R/(x_1, \ldots, x_d)$ is regular,
so is $R_{\p}/(x_1, \ldots, x_d)_{\p}$,
whence $x_1, \ldots, x_d$ is part of a (usual) regular system of
parameters in $R_{\sp}$.

\lemma
\label{\lmlbvjac}
Let $R$ be a regular domain
and $x_1, \ldots, x_d$ a generalized regular system of parameters.
Then for any normalized valuation $v$ as in the definition of adjoints,
$v(J_{R_v/R}) \ge v(x_1 \cdots x_d) - 1$.
\endb

\proof
By possibly taking a subset of the $x_i$,
without loss of generality all $v(x_i)$ are positive.
Let $\p$ be the contraction of the maximal ideal of $R_v$ to $R$.
After localizing at $\p$,
$x_1, x_2, \ldots, x_d$ are a part of a regular system of parameters
(see comment above the lemma).
We may possibly extend the $x_i$ to a full regular system of parameters,
so we may assume that $\p = (x_1, \ldots, x_d)$
is the unique maximal ideal in $R$.
We may also assume that $v(x_1) \ge v(x_2) \ge \cdots \ge v(x_d) \ge 1$.

If $d = 0$,
the lemma holds trivially.
If $d = 1$,
then $v$ is the $\p$-adic valuation,
in which case $R_v = R$,
$J_{R_v/R} = R$.
As $v$ is normalized,
$v(x_1) = 1$,
and the lemma holds again.
Now let $d > 1$ and let $S = R[{x_1 \over x_d}, \ldots, {x_{d-1} \over x_d}]$.
Then $S$ is a regular ring contained in $R_v$,
and ${x_1 \over x_d}, \ldots, {x_{d-1} \over x_d}, x_d$
are a generalized regular system of parameters in $S$.
Clearly $J_{S/R}$ equals $x_d^{d-1}S$.
By induction on $\sum_i v(x_i)$
we conclude that $v(J_{R_v/S}) \ge
v({x_1 \over x_d} \cdots {x_{d-1} \over x_d} x_d) - 1$,
so that
by Lipman and Sathaye~\cite[page 201]{LipSat},
$v(J_{R_v/R}) = v(J_{S/R}) + v(J_{R_v/S})
\ge v(x_d^{d-1}) + v({x_1 \over x_d} \cdots {x_{d-1} \over x_d} x_d) - 1
= v(x_1 \cdots x_d) - 1$.
\qed

Though in general not as nicely behaved as variables in a polynomial
or power series ring, generalized regular systems of parameters come
close to them in many aspects.
One interesting property is the following.

\prop
\label{\lemmonomrep} Let $R$ be a regular domain and let $x_1, \dots, x_d$ be a 
generalized regular system of parameters of $R$.
Furthermore let $s \leq d$, $\p = (x_1, \dots, x_s)$ and
let $f$ be a non-zero element of $R$.
Then there exist monomials $m_1, \dots, m_t$ in $x_1, \dots, x_s$ and elements
$h, g_1, \dots, g_t \in R \setminus \p$ such that
  $$ h \cdot f = \sum_{i=1}^t g_i \cdot m_i. $$
\endb

\proof
Clearly we may assume that $R$ is local with maximal ideal $\p$,
and we then prove the proposition with $h = 1$.
First we reduce to the case of complete local rings:
Let $\widehat{R}$ be the completion of $R$, and note that $x_1,
\dots, x_s$ is a regular system of parameters of $\widehat{R}$.
Suppose we know the result for $\widehat{R}$ and $x_1, \dots,
x_s \in \widehat{R}$.
Write
  $ f = \sum_{i = 1}^t h_i m_i  $
for some $h_i \in \widehat{R}$, $h_i \notin \p$.
Clearly we may assume that none of the monomials is a multiple of another one.
Let $I = (m_1, \dots, m_t) \subseteq R$.
As
  $ f \in I \widehat{R} \cap R = I $
by faithful flatness, we may write
  $ f = \sum_{i=1}^t g_i m_i $
with $g_i \in R$, and in $\widehat{R}$ we get
  $ \sum_{i=1}^t (g_i - h_i) m_i = 0 $
hence we conclude from \cite{Kap}, \S 5, that $g_i - h_i \in \p 
\widehat{R}$, implying that $g_i \notin \p$.
Thus it suffices to prove the proposition
in the case $R$ is complete local with maximal ideal
$\p = (x_1, \ldots, x_s)$.

Assume now that $R$ is complete and let $f \in R$.
Assume $f \in \p^{n_1} \setminus \p^{n_1+1}$.
Then we may write
  $$ f = \sum_{i=1}^{t_1} a_{1i} m_{1i} + f_{2} $$
with some (unique) monomials $m_{1i}$ of degree $n_1$ in $x_1, \dots, x_s$ and some 
$a_{1i} \notin \p$ (unique mod $\p$)
and with some $f_2 \in \p^{n_1+1}$.
Let $M_1 = (m_{11}, \dots, m_{1t_1})$.
If $f_2 = 0$ we are done, otherwise we write
  $$ f_2 = \sum_{i=1}^{t_2} a_{2i} m_{2i} + f_3 $$
with some $a_{2i} \notin \p$,
some monomials $m_{2i}$ of degree $n_2$ in $x_1, \dots, x_s$,
and some $f_3 \in \p^{n_2+1}$.
Set $M_2 = M_1 + (m_{21}, \dots, m_{2t_2})$
and continue.
In this way we get an ascending chain $M_1 \subseteq M_2 \subseteq \cdots $ of 
monomial ideals, which must stabilize eventually,
$M_{\rho} = M_{\rho + 1} = \cdots =: M_{\infty}$.
Let
  $$M_{\infty} = (m_1, \dots, m_t), $$
with each $m_i$ a monomial of degree $d_i$ in $x_1, \dots, x_s$.
We may assume that none of these monomials divides any of the other ones
and that all $m_i$ appear in a presentation of some $f_j$ as above.
Then in each step above we may write
  $$ f_l = \sum_{i = 1}^t n_{li} m_i + f_{l+1} $$
with $n_{li} \in \p^{n_l - d_i}$,
and where furthermore if $l$ is the smallest integer such that
$m_i$ appears with a non-trivial coefficient in the expansion of $f_l$,
we have $n_{li} \not \in \p$.
Hence
  $$ f = \sum c_{li} m_i + f_{l+1} $$
with some $c_{li} \notin \p$ (or $c_{li} = 0$), and with 
$c_{l+1,i} - c_{li} \in \p^{n_l - d_i}$ (and $f_{l+1} \in \p^{n_{l+1}+1}$).
As $R$ is complete, this converges, and we get
  $$ f = \sum c_i m_i $$
with some $c_i \notin \p$.
\qed

\section{Integral closures of (general) monomial ideals}

Monomial ideals typically denote
ideals in a polynomial ring or in a power series ring over a field
that are generated by monomials in the variables.
Such ideals have many good properties,
and in particular,
their integral closures and multiplier ideals are known to be monomial as well.
The just stated result on multiplier ideals for the standard monomial ideals is due to Howald~\cite{Howald}.
In this section we consider generalized monomial ideals and present 
their integral closures.
For alternate proofs on the integral closure of generalized monomial ideals see Kiyek and St\"uckrad~\cite{KiSt}.

We define monomial ideals more generally:

\defin
Let $R$ be a regular domain,
and let $x_1, \ldots, x_d$ in $R$ be a generalized regular system of parameters.
By a {\bf monomial ideal (in $\bfi x_{\bfis 1}, \bfi{$\ldots$}, \bfi x_{\bfis d}$)}
we mean an ideal in $R$ generated by monomials in $x_1, \ldots, x_d$.
\endb

As in the usual monomial ideal case,
we can define the Newton polyhedron:

\defin
Let $R, x_1, \ldots, x_d$ be as above,
and let $I$ be an ideal generated by monomials
$\underline x^{\underline a_1}, \ldots, \underline x^{\underline a_s}$.
Then the {\bf Newton polyhedron of $I$}
(relative to $x_1, \ldots, x_d$) is the set
$$
\NP(I) = \{\underline e \in \bbQ_{\geq 0}^d \,|\,
\underline e \ge \Sigma_i c_i \underline a_i, \text{ for some }
c_i \in \bbQ_{\ge 0}, \Sigma_i c_i = 1\}.
$$
$\NP(I)$ is an unbounded closed convex set in $\bbQ_{\ge 0}^d$.
We denote the interior of $\NP(I)$ as $\NPint (I)$.
\endb

\thm
\label{\thmintclmonom}
Let $R$ be a regular domain
and $x_1, \ldots, x_d$ a generalized regular system of parameters.
Let $I$ be an ideal generated by monomials
in $x_1, \ldots, x_d$.
The integral closure $\ic {I^n}$ of $I^n$ equals
$$
\ic {I^n} = ( \{ \underline x^{\underline e} \,|\, \underline e \in n \cdot 
\NP(I) \cap \bbN^d \} ),
$$
so it is generated by monomials.
\endb

\proof
As $\NP(I^n) = n \cdot \NP(I)$, we may assume that $n = 1$. Write
$I = (\underline x^{\underline a_1}, \ldots, \underline x^{\underline a_s})$.
Let $\alpha = \underline x^{\underline e}$ be
such that $\underline e \in \NP(I) \cap \bbN^d$.
Then there exist $c_1, \ldots, c_s \in \bbQ_{\ge 0}$
such that $\sum c_i = 1$ and $\underline e \ge \sum c_i \underline a_i$
(componentwise).
Write $c_i = m_i/n$ for some $m_i \in \bbN$ and $n \in \bbN_{>0}$.
Then
$$
\alpha^n =
x_1^{n e_1 - \Sigma m_i a_{i1}}
\cdots
x_d^{n e_d - \Sigma m_i a_{id}}
(\underline x^{\underline a_1})^{m_1}
\cdots
(\underline x^{\underline a_s})^{m_s} \in I^{m_1+ \cdots + m_s} = I^n,
$$
so that $\alpha \in \ic I$.
It remains to prove the other inclusion.

Let $S$ be the set of hyperplanes that bound $\NP(I)$
and are not coordinate hyperplanes.
For each $H \in S$,
if an equation for $H$ is $h_1 X_1 + \cdots + h_d X_d = h$
with $h_i \in \bbN$ and $h \in \bbN_{>0}$,
define
$I_H = (\underline x^{\underline e} ~|~
\underline e \in \bbN^d, \sum_i h_i e_i \ge h)$.
Clearly $I \inc I_H$,
$\NP(I_H) \subseteq \{ e \in \bbQ_{\geq 0}^d: \sum h_i e_i \geq h \}$,
and
$\NP(I_H) \cap \bbN^d =
\{ e \in \bbQ_{\geq 0}^d: \sum h_i e_i \geq h \} \cap \bbN^d$.
Suppose that the theorem is known for the (generalized) monomial ideals $I_H$.
Then
$$
\eqalignno{
\ic I &\inc \bigcap_{H \in S} \ic{I_H} \cr
&\inc \bigcap_{H \in S} \left( \{ \underline x^{\underline e} ~|~
\underline e \in \bbN, \hbox{$\sum_i h_i e_i \ge h$},
\hbox{ if $H$ is defined by $\sum_i h_i X_i = h$} \} \right) \cr
&= ( \{ \underline x^{\underline e} \,|\, \underline e \in \NP(I) 
\cap \bbN^d \} ). \cr
}
$$
Thus it suffices to prove the theorem for $I_H$.
As before,
let $\sum_i h_i X_i = h$ define $H$.
By possibly reindexing,
we may assume that $h, h_1, \ldots, h_t$ are positive integers
and that $h_{t+1} = \cdots = h_d = 0$.
As noted above it suffices to show
  $$ \overline{I_H} = ( \{ \underline x^{\underline e} \,|\, e \in \bbN^d
     \text{ and } \sum h_i e_i \geq h \}  ). $$

Let $Y_1, \ldots, Y_t$ be variables over $R$
and $R' = R[Y_1, \ldots, Y_t]/(Y_1^{h_1} - x_1, \ldots, Y_t^{h_t} - x_t)$.
This is a free finitely generated $R$-module
and $Y_1, \ldots, Y_t$ is a regular sequence in $R'$.
Set $\p = (Y_1, \ldots, Y_t)R'$.
Then $R'/\p = R/(x_1, \ldots, x_t)$ is a regular domain,
so $\p$ is a prime ideal,
and for any prime ideal $\q$ in $R'$ containing $\p$,
$R'_{\sq}$ is a regular local ring.
By construction,
$I_H R'$ is contained in $(Y_1, \ldots, Y_t)^h = \p^h$.
As $R'_{\sp}$ is a regular local ring,
$\p^h R'_{\sp}$ is integrally closed,
and as $R'$ is finitely generated over a locally formally equidimensional
(regular) ring,
$R'_{\sq}$ is locally formally equidimensional
for every prime ideal $\q$ containing $\p$.
By a theorem of Ratliff, from~\cite{Ra74},
since $\p$ is generated by a regular sequence,
the integral closure of $\p^h R'_{\sq}$ has no embedded prime ideals.
It follows that the integral closure of $\p^h R'_{\sq}$ is $\p^h R'_{\sp}
 \cap R'_{\sq}$.
As $R'_{\sq}$ is a regular domain
and $\p$ is generated by a regular sequence,
$\p^h R'_{\sp} \cap R'_{\sq} = \p^h R'_{\sq}$.
It follows that $\p^h R'_{\sp} \cap R' = \p^h$ is the integral closure of $\p^h$.
Hence $\ic {I_H} \inc \ic{\p^h} \cap R = \p^h \cap R$,
and by freeness of $R'$ over $R$,
the last ideal is exactly
$(\underline x^{\underline e} \,|\,
\underline e \in \bbN, \sum_i h_i e_i \ge h)$,
which finishes the proof.
\qed

\section{Rees valuations of (general) monomial ideals}

Recall that the Rees valuations of a non-zero ideal in a Noetherian domain
form a unique minimal set $\Reesval(I)$ of finitely many normalized valuations
such that for all positive integers $n$,
$\ic {I^n}
= \{r \in R ~|~ v(r) \ge n v(I) \hbox{ for all } v \in \Reesval(I)\}$.

In an arbitrary Noetherian domain,
for arbitrary ideals $I$ and $J$,
$\Reesval(I) \cup \Reesval(J) \inc \Reesval(IJ)$,
and equality holds in two-dimensional regular domains.
(This has appeared in the literature in several places,
see for example~Muhly-Sakuma~\cite{MS63},
or the Rees valuations chapter in the upcoming book~\cite{SwHu}.)

We will prove that the Rees valuations of an ideal
generated by monomials in a regular system of parameters
are especially nice.

\defin
Let $R$ be a regular domain,
and let $x_1, \ldots, x_d$ be a generalized regular system of parameters.
A valuation $v$ on the field of fractions of $R$
is said to be {\df monomial} on $x_1, \ldots, x_d$
if for some $i_1, \ldots, i_s \in \{1, \ldots, d\}$,
for any polynomial
$f = \sum c_{\nu} x_{i_1}^{\nu_{i_1}} \cdots x_{i_s}^{\nu_{i_s}} \in R$
with all $c_{\nu}$ either $0$ or not in $(x_{i_1}, \ldots, x_{i_s})$,
$v(f) = \min\{v({\underline x}^{\underline \nu}) ~|~ c_{\nu} \not = 0\}$.
When the $x_i$ are understood from the context,
we say that $v$ is {\df monomial}.
\endb

Observe that $v(f) = 0$ for any $f \notin (x_{i_1}, \dots, x_{i_s})$.
In particular $v(x_j) = v(1) = 0$ if $j \not \in \{i_1, \ldots, i_s\}$.

\prop
\label{\propmonomialval} Let $R$ be a regular domain, let $x_1, \dots ,
x_d$ be a generalized regular system of parameters, and let $a_1,
\dots, a_d$ be non-negative rational numbers, not all of them zero.
Then there exists a unique valuation $v$ on the field of fractions $K = Q(R)$ of 
$R$
that is monomial on $x_1, \dots, x_d$, with $v(x_i) = a_i$.
\endb

\proof By reindexing we may assume that
  $a_1 > 0, \ldots, a_s > 0, \quad a_{s+1} = \cdots = a_d = 0$
for some $s > 0$,
and we also may assume that all $a_i$ are integers.

The uniqueness of $v$ is immediate by Proposition~\lemmonomrep.
To prove the existence we may replace $R$ by $R_{\sp}$ (with $\p = (x_1, 
\ldots, x_s)$) and assume that $R$ is local.
Let $R'$ be the regular local ring obtained by adjoining a $a_i^{\text{\sevenrm 
th}}$-root $y_i$ of $x_i$ to $R$ ($i = 1, \dots, s$) and
let $\n$ be the maximal ideal of $R'$.
Then the $\n$-adic valuation $w$ on
$L = Q(R')$ is monomial in $y_1, \dots, y_s$ with $w(y_i) = 1$ for all
$i$.
%
The restriction $v := w \vert_K$ is a monomial valuation as desired.
\qed

\cor
\label{\corcorr}
Let $R$ be a regular domain,
and let $x_1, \ldots, x_d$ be a generalized regular system of parameters.
Let $I$ be an ideal generated by monomials in $x_1, \ldots, x_d$.
Then all the Rees valuations of $I$ are monomial in $x_1, \ldots, x_d$.
Furthermore,
if $H_1, \dots , H_{\rho}$ are the non-coordinate hyperplanes bounding $\NP(I)$,
then the $H_j$ are in one-to-one correspondence
with the Rees valuations $v_j$ of $I$.
\endb

\proof
The Newton polyhedron $\NP(I)$ of $I$
is the intersection of finitely many half-spaces in $\bbQ^d$. 
Some of them are coordinate half-spaces $\{x_i \geq 0\}$, 
each of the others is determined by a
hyperplane $H$ of the form $h_1 x_1 + \cdots + h_d x_d = h$,
with $h_1, \ldots, h_d, h$ non-negative integers, $h > 0$,
and $\gcd(h_1, \ldots, h_d, h) = 1$.
This hyperplane corresponds to a valuation $v_H$
that is monomial on $x_1, \ldots, x_d$
and such that $v_H(x_i) = h_i$.
By Theorem~{\thmintclmonom},
the integral closure of $I$ is determined by these $v_H$.
Using $\NP(I^n) = n \cdot \NP(I)$,
we see that the integral closure of $I^n$ is also determined by these $v_H$.
So each Rees valuation is one such $v_H$.
Suppose that the set of Rees valuations is a proper subset of the set
of all the $v_H$.
Say one such $v_H$ is not needed in the computation of the
integral closures of powers of $I$.
Since the hyperplanes $H$ were chosen to be irredundant,
by omitting any one of them,
we get a point $(e_1, \ldots, e_d) \in \bbQ_{\ge 0}^d$
which is on the unbounded side of all the hyperplanes bounding $\NP(I)$
other than $H$, but is not on the unbounded side of $H$.
There exists $m_1, \ldots, m_d \in \bbN$ and $n > 0$
such that for each $i$,
$e_i = m_i/n$.
Then by assumption,
$x_1^{m_1} \cdots x_d^{m_d} \in \ic{I^n}$, but $(m_1, \ldots, m_n) 
\notin n \cdot \NP(I)$, a contradiction.
\qed

The following is a local version of Howald~\cite[Lemma 1]{Howald}.
Howald's proof relies on the existence of a log resolution.

\lemma
\label{\lmadjmonom}
Let $R$ be a regular domain,
and let $x_1, \ldots, x_d$ be a generalized regular system of parameters.
Let $v$ be a discrete valuation
that is monomial on $x_1, \ldots, x_d$, non-negative on $R$,
and has value group contained in $\bbZ$.
Then
$$
v(J_{R_v}/R) = v(x_1\cdots x_d) - gcd(v(x_i) | i).
$$
\endb

\proof
Since $v$ is monomial in the $x_i$,
the center of $v$ on $R$
is contained in $\m = (x_1, \ldots, x_d)$.
By localizing,
we may assume that $\m$ is the only maximal ideal in $R$.
Let $a_i = v(x_i)$.
Without loss of generality $a_1 \ge a_2 \ge \cdots \ge a_d$,
and let $s$ be the largest integer such that $a_s > 0$.
As $v$ is monomial, if $s = 0$,
then $v = 0$ and the lemma holds trivially.
Thus we may assume that $s > 0$.
If $s = 1$, necessarily $a_1 = \gcd(v(x_i) | i)$,
and $v$ is $a_1$ times the $(x_1)$-grading.
Then $R_v = R_{(x_1)}$, $J_{R_v/R} = R_v$,
$v (J_{R_v/R}) = 0 = v(x_1 \cdots x_d) - \gcd(v(x_i) | i)$.
So the lemma holds in the case $s = 1$.
We proceed by induction on $\sum_i a_i$.
We may assume that $s > 1$.
Let $S = R[{x_1 \over x_s}, \ldots, {x_{s-1} \over x_s}]$.
Then $S$ is a regular ring contained in $R_v$,
${x_1 \over x_s}, \ldots, {x_{s-1} \over x_s}, x_s, \ldots, x_d$
is a generalized regular system of parameters.
For these elements in $S$,
$v$ is still a monomial valuation,
their $v$-values are non-negative integers,
and the total sum of their $v$-values is strictly smaller than $\sum_i a_i$.
Thus by induction,
$$
\eqalignno{
v(J_{R_v/S})
&= v\left({x_1 \over x_s} \cdots {x_{s-1} \over x_s} x_s \cdots x_d\right)
- \gcd(v(x_1/x_s), \ldots, v(x_{s-1}/x_s), v(x_s), \ldots, v(x_d)) \cr
&= v(x_1 \cdots x_d) - (s-1) v(x_s)
- \gcd(v(x_1), \ldots, v(x_d)). \cr
}
$$
As $R \inc S \inc R_v$ are all finitely generated algebras over $R$
that are regular rings and have the same field of fractions,
by Lipman and Sathaye~\cite[page 201]{LipSat},
$J_{R_v/R} = J_{S/R} J_{R_v/S}$.
Clearly $J_{S/R}$ equals $x_s^{s-1}$,
whence
$$
\eqalignno{
v(J_{R_v/R})
&= v(x_s^{s-1}) +
v(x_1 \cdots x_d) - (s-1) v(x_s) - \gcd(v(x_1), \ldots, v(x_d)) \cr
&= v(x_1 \cdots x_m) - \gcd(v(x_1), \ldots, v(x_d)).
& \qedbox \cr
}
$$

\section[\sectgenmonRVa]{Adjoints of (general) monomial ideals}

A proof similar to the proof of Theorem~{\thmintclmonom}
shows that the adjoint of a (general) monomial ideal is monomial.
This generalizes Howald's result~\cite{Howald}.

\thm
\label{\thmadjmonomval}
Let $R$ be a regular domain,
and let $x_1, \ldots, x_d$ be a generalized regular system of parameters.
Let $I$ be an ideal generated by monomials in $x_1, \ldots, x_d$.
Then for all $n \ge 1$,
$$
\eqalignno{
\adj (I^n) &= \bigcap_v\, (\{\underline x^{\underline e} ~|~
v(\underline x^{\underline e}) \ge v(I^n) - v(x_1\cdots x_d) + 1\}) \cr
&= \bigcap_v\, (\{\underline x^{\underline e} ~|~
v(\underline x^{\underline e}) \ge v(I^n) - v(J_{R_v/R})\}) \cr
&= ( \{ \underline x^{\underline e}: \underline e \in \bbN^d \text{ and }
e + (1, \ldots 1) \in \NPint(I^n) \} ), \cr
}
$$
as $v$ varies over the (normalized) Rees valuations of $I$.
In particular,
the adjoint is also generated by monomials.
\endb

\proof
As $I^n$ is monomial
and as the Rees valuations of $I^n$
are contained in the set of Rees valuations of $I$,
it suffices to prove the theorem for $n = 1$.
By Corollary~{\corcorr} and by Lemma~\lmadjmonom,
the second and third equalities hold.
So it suffices to prove that $\adj I$ equals the other three expressions (when $n = 1$).

First we prove that $\underline x^{\underline e} \in \adj(I)$
whenever $\underline e \in \bbN^d$
with $\underline e + (1,\ldots, 1) \in \NPint(I))$.
Let $v$ be a valuation as in the definition of $\adj(I)$.
As $(x_1 \cdots x_d \underline x^{\underline e})^n \in \ic {I^{n+1}}$
for some positive integer $n$,
$v(x_1 \cdots x_d \underline x^{\underline e}) > v(I)$.
As $v$ is normalized,
$v(\underline x^{\underline e}) \ge v(I) - v(x_1 \cdots x_d) + 1$.
By Lemma~{\lmlbvjac},
$v(J_{R_v/R}) \ge v(x_1 \cdots x_d) - 1$,
so that $v(\underline x^{\underline e}) \ge v(I) - v(J_{R_v/R})$.
As $v$ was arbitrary,
this proves that
$(\underline x^{\underline e} \,|\,
\underline e \in \bbN^d,
\underline e + (1,\ldots, 1) \in \NPint(I)) \subset \adj I$.
It remains to prove the other inclusion.

Let $S$ be the set of bounding hyperplanes of $\NP(I)$
that are not coordinate hyperplanes.
For each $H \in S$,
if an equation for $H$ is $h_1 X_1 + \cdots + h_d X_d = h$
with $h_i \in \bbN$ and $h \in \bbN_{>0}$,
define
$I_H = (\underline x^{\underline e} ~|~
\underline e \in \bbN^d, \sum_i h_i e_i \ge h)$.
By the definition of Newton polyhedrons,
$I \inc I_H$.

By possibly reindexing,
without loss of generality $h_1, \ldots, h_t > 0$
and $h_{t+1} = \cdots = h_d = 0$.
By Proposition~{\propmonomialval}
there exists a monomial valuation $v_H$ on $Q(R)$
defined by $v_H(x_i) = h_i$.
By construction,
$v_H(I) \ge v_H(I_H) \ge h$ (even equalities hold),
and
$\adj(I_H) \inc \{r \in R ~|~ v_H(r) \ge v_H(I_H) - v_H(J_{R_{v_H}/R})\}$.
By the properties of $v_H$,
the last ideal is generated by monomials in the $x_i$.
By Lemma~\lmadjmonom,
$v_H(J_{R_{v_H}/R}) = v_H(x_1 \cdots x_d) - 1$,
so that
$$
\eqalignno{
\adj(I_H)
&\inc ( \{ \underline x^{\underline e} ~|~
\underline e \in \bbN^d,
v_H(\underline x^{\underline e}) > v_H(I_H) - v_H(x_1 \cdots x_d) \} ) \cr
&\inc ( \{ \underline x^{\underline e} ~|~
\underline e \in \bbN^d,
\Sigma_i h_i (e_i + 1) > v_H(I_H) \} ) \cr
&\inc ( \{ \underline x^{\underline e} ~|~
\underline e \in \bbN^d,
\Sigma_i h_i (e_i + 1) > h \} ), \cr
}
$$
whence
$$
\eqalignno{
\adj I &\inc \bigcap_{H \in S} \adj(I_H) \cr
&\inc \bigcap_{H \in S} \left( \{ \underline x^{\underline e} ~|~
\underline e \in \bbN,
\Sigma_i h_i (e_i + 1) > h,
\hbox{ if $H$ is defined by $\sum_i h_i X_i = h$}\right \} ) \cr
&= \left( \{ \underline x^{\underline e} ~|~
\underline e \in \bbN^d,
\underline e + (1, \ldots, 1) \in \NPint(I) \} \right). & \qedbox \cr
}
$$

This theorem allows to address the subadditivity problem for monomial ideals:

\cor
\label{\corsubaddmon}
 Let $I, J \subseteq R$ be ideals generated by monomials in
the generalized regular system of parameters $x_1, \ldots , x_d$.
Then
  $$ \adj (I  J) \subseteq \adj(I) \cdot \adj(J). $$
\endb

\proof Let $\underline x^{\underline a} \in \adj(I  J)$ be a monomial.
By Theorem~\thmadjmonomval,
$\underline a + (1, \ldots, 1) \in \NPint(I \cdot J)$.
As $\NPint(I  J) \subseteq \NPint(I) + \NPint(J)$, there exist 
$\underline b \in \NPint(I)$ and $\underline c \in \NPint(J)$ with 
$\underline a + (1, \ldots, 1) = \underline b + \underline c$.
Set $\underline f = (f_1, \dots, f_d)$
and $\underline g = (g_1, \dots, g_d)$ with
$f_i = \lceil b_i \rceil - 1$ 
and $g_i = \lfloor c_i \rfloor$.
Then $\underline x^{\underline g}$ and $\underline x^{\underline f}$
are monomials with $\underline x^{\underline g} \cdot
\underline x^{\underline f} = \underline x^{\underline a}$, and
furthermore
  $$ \eqalignno{
  \underline f + (1, \ldots, 1) \, & \, \in \underline b +
  \bbQ^d_{\geq 0} \subseteq \NPint(I), \cr
  \underline g + (1, \ldots, 1) \, & \, \in \underline c +
  \bbQ^d_{\geq 0} \subseteq \NPint(J),
  } $$
implying by Theorem~{\thmadjmonomval}
that $\underline x^{\underline f} \in \adj(I)$ and
$\underline x^{\underline g} \in \adj(J)$.
\qed

From the proof of Theorem~\thmadjmonomval,
it is clear that the Rees valuations of the adjoint depend on the
Rees valuations of the original ideal.
The number of Rees valuations of $I$
need not be an upper bound on the number of Rees valuations of $\adj(I)$,
and there is in general no overlap between the set of Rees valuations of $I$
and the set of Rees valuations of $\adj I$.

\example
Let $R$ be a regular local ring with regular system of parameters $x, y$.
Let $I$ be the integral closure of $(x^5, y^7)$.
Then by the structure theorem,
$I$ has only one Rees valuation,
and $I = (x^5, x^4y^2, x^3y^3, x^2y^5, xy^6, y^7)$.
By \cite{HunSwa}, by~\cite{Howald}, or by~Theorem~\thmadjmonomval,
$\adj(I) = (x^4, x^3y, x^2y^2, xy^4, y^5)$,
which is not the integral closure of $(x^4, y^5)$.
Thus $\adj(I)$ has more than one Rees valuation.
In fact,
it has two Rees valuations, both of which are monomial
and neither of which is equivalent to the Rees valuation of $I$:
$v_1(x) = 1 = v_1(y), v_1(\adj(I)) = 4$,
and $v_2(x) = 3, v_2(y) = 2, v_2(\adj(I)) = 10$.

Nevertheless,
the one Rees valuation of $I$ still determines the adjoints of all the
powers of $I$.

\section[\sectchar]{Adjoints of ideals and Rees valuations}

In this section we characterize those ideals $I$ for which $\adj(I^n)$ is 
determined by the Rees valuations of $I$ for all $n$. In the last section we 
have seen that this is the case for monomial ideals. 
That the Rees valuations of an ideal $I$ should
play a crucial role in determining the adjoint of $I$ in general is also implied 
by the following result:

\prop
\label{\propreesvalcon}
Let $I$ be an ideal in a regular domain $R$,
and let $\val$ be a finite set of valuations on the field of fractions of $R$
such that for all $n \in \bbN$,
$$
\adj I = \bigcap_{v \in \sval} \, \{r \in R ~|~ v(r) \ge v(I) - v(J_{R_v/R}) \}.
$$
Then $\val$ contains the Rees valuations of $I$.
\endb

\proof
Assume that there exist some Rees valuations of $I$ not contained in $\val$. By the defining property of Rees valuations there exist a non-negative integer $n$
and an element $r \in R$ with

\item{(1)}  $v(r) \geq n \cdot v(I) \, $ for all $ \, v \in \val$.
\item{(2)}  $r \notin \ic{I^n}$.

Let $w$ be a Rees valuation of $I$ with $w(r) \leq n \cdot w(I) - 1$. Assume that $I$ is $l$-generated and let $t \geq l \cdot w(I)$. Then
  $$ w(r^t) = t \cdot w(r) < (nt - l + 1) w(I),$$
hence
  $$ r^t \notin \ic{I^{nt-l+1}}. $$
On the other hand,
  $$ v(r^t) \geq nt \cdot v(I) \geq nt \cdot v(I) - v( J_{R_v/R}) \, \,
    \text{ for all } \, v \in \val, $$
implying that 
  $$ r^t \in \adj(I^{nt}) \subseteq \ic{I^{nt-l+1}} $$
by~\cite{Lipmanadj}, (1.4.1), a contradiction.
\qed

It is not true in general that the set of Rees valuations
determines the adjoint of an arbitrary ideal:

\example
Let $(R,\m)$ be a $d$-dimensional regular local ring,
with $d > 2$,
and let $\p$ be a prime ideal in $R$ of height $h \in \{2, \ldots, d-1\}$
generated by a regular sequence.
Then the $\p$-adic valuation $v_{\sp}$ is the only Rees valuation of $\p$.
If $v_{\sp}$ defined $\adj(\p^n)$ in the sense that
  $$ \adj(\p^n) = \{r \in R ~|~ v_{\sp}(r) \ge n v_{\sp}(\p) - 
     v_{\sp}(J_{R_{v_{\ssp}}/R})\} \quad  \text{ for all } n,$$
then, as $v_{\sp}(\p) = 1$ and $J_{R_{v_{\ssp}}/R} = \p^{h-1} R_{v_{\ssp}}$,
it follows that
  $$ \adj(\p^{h-1}) = \{r \in R ~|~ v_{\sp}(r) \ge 0\} = R. $$
However, if $\p$ is generated by elements in $\m^{e}$,
where $e \ge d/(h-1)$,
and if $v$ denotes the $\m$-adic valuation,
then
  $$ \adj(\p^{h-1}) \inc \{r \in R ~|~ v(r) \ge v(\p^{h-1}) - 
v(J_{R_v/R})\} \inc \{r \in R ~|~ v(r) \ge d - (d-1)\} \inc \m, $$
which is a contradiction.
A concrete example of this is $R = k[[X,Y,Z]]$
with the prime ideal $\p = (X^4 - Z^3, Y^3 - X^2 Z)$,
which defines the monomial curve $(t^9, t^{10}, t^{12})$.

\vskip 3pt

The following is a geometric reformulation of~\cite{Rees},
see also~\cite{CHRR}, 2.3 or~\cite{Lipproxy}, 1.4:

\remark
\label{\rmkreesval}
Let $R$ be a regular domain
and let $I \subseteq R$ be an ideal of $R$.
Let $Y = \text{Spec}(R)$,
$P = R[IT]$,
the Rees ring of $I$ and let
$\overline{P}$ be its normalization and
$\varphi: X = \text{Proj}(\overline{P}) \rightarrow Y$ the induced scheme. Then $X/Y$ is essentially of finite type by~\cite{LipSat}, p.~200 (see also~\cite[9.2.3]{SwHu}, for details). Thus $\varphi$ is a projective, 
birational morphism, $X$ is a normal, Noetherian scheme and $I \hbox{\Cal O}_X$ is an 
invertible ideal.
Let
$\P_1, \dots , \P_r$ be the irreducible components of $\V(I)$ (i.e., those
points $x$ of $X$ of codimension 1 such that $I \hbox{\Cal O}_{X,x}$ is a proper
ideal of $\hbox{\Cal O}_{X,x}$).
Then $\hbox{\Cal O}_{X, \sP_i}$ is a discrete valuation
ring (with field of fractions $K = Q(R)$) and the corresponding valuations
$v_1, \dots, v_r$ are exactly the Rees valuations of $I$.

If $(R, \m)$ is local and $I$ is $\m$-primary,
the Rees valuations correspond to the irreducible
components of the closed fibre $\varphi^{-1}(\m)$
which in this case is a $(\text{dim}(R)-1)$-dimensional projective
scheme (in general neither reduced nor irreducible).

\bigskip

Let $f : Z \rightarrow Y$ be birational and of finite type.
Then the Jacobian ideal $\hbox{\Cal J}_{Z/Y} \subseteq \hbox{\Cal O}_Z$ is well-defined
(being locally the $0$\th-Fitting ideal of the relative K\"ahler differentials).
If in addition $Z$ is normal, then
$$
 \omega_{Z/Y} := \hbox{\Cal O}_Z : \hbox{\Cal J}_{Z/Y} 
 = \hbox{\Cal H}om_Z ( \hbox{\Cal J}_{Z/Y}, 
 \hbox{\Cal O}_Z)
$$
is a canonical dualizing sheaf for $f$ with
$$
\hbox{\Cal O}_Z \subseteq \omega_{Z/Y} \subseteq \hbox{\Cal M}_Z,
$$
where $\hbox{\Cal M}_Z$ denotes the constant sheaf of meromorphic functions on $Z$.
If
$$
g: Z' \longrightarrow Z
$$
is another birational morphism and if $g$ is proper and $Z'$ is normal as well,
then
$$
g_* \omega_{Z'/Y} \subseteq \omega_{Z/Y}
$$
(cf.~\cite[2.3]{LipSat} and~\cite[\S 4]{LipTei}).

\thm
\label{\thmadjone}
Let $R$ be a regular domain and let $I \subseteq R$
be a non-trivial ideal. Furthermore let $Y = \text{Spec}(R)$ and
$\varphi: X \rightarrow Y$  be the normalized blow-up of $I$.
Then the following are equivalent:
\item{(1)}
$\adj(I^n)
= \displaystyle\bigcap_{v \in \sReesval(I)} \, \{ r \in R ~|~ v(r) \ge n \cdot v(I) 
- v(J_{R_v/R})\}$
for all positive integers $n$.
\item{(2)}
If $Z$ is a normal scheme and $\pi: Z \rightarrow X$ is proper and birational,
then
$$
\pi_* \omega_{Z/Y} = \omega_{X/Y}.
$$
\endb

\remark
\label{\rmkpseudorat}
If in the situation of~\thmadjone (2) the scheme $X$ is Cohen--Macaulay as well,
then $X$ has pseudo-rational singularities only
(\cite[\S 4]{LipTei}).

\remark
\label{\rmkreesvaluni}
In the situation of~\thmadjone (1) the set $\Reesval(I)$ is the unique smallest set of valuations defining $\adj(I^n)$ in view of Proposition~\propreesvalcon.

\medskip

\proofof{Theorem~\thmadjone}
If $f: Z\rightarrow Y$ is proper and birational with $Z$ normal
and $I \hbox{\Cal O}_Z$ invertible,
we set
$$
\text{adj}_Z(I^n) = H^0(Z, I^n \omega_{Z/Y}) \quad (\, \subseteq R  \, ).
$$
Then $\adj(I^n) = \bigcap \text{adj}_Z(I^n)$ by~\cite{Lipmanadj}, where $f: Z\rightarrow Y$ varies over all such morphisms.
By the universal properties of blow-up and normalization, any such
$f$ factors as
$$
Z\overarrow{\pi} X \overarrow{\varphi} Y,
$$
As $ \pi_* \omega_{Z/Y} \inc \omega_{X/Y},$
and as $I \hbox{\Cal O}_Z$ is invertible,
this implies by the projection formula

$$
\eqalignno{
   H^0(Z, I^n \omega_{Z/Y}) \, \, & = H^0(X, \pi_*   I^n \omega_{Z/Y} ) \cr
   & = H^0(X, I^n \pi_* \omega_{Z/Y}) \cr
   & \inc H^0(X, I^n \omega_{X/Y}), \cr
}
$$
and therefore 

\vskip 2pt
\item{($*$)} $\qquad \qquad \qquad \adj_Z(I^n) \subseteq \text{adj}_X(I^n) \quad 
         \text{ for all positive integers } n  $ 
\vskip 2pt

\noindent for any such $f: Z \rightarrow Y$.

As $\omega_{X/Y}$ is reflexive by~\cite[p. 203]{LipSat},
and as $I \hbox{\Cal O}_X$ is invertible,
$I^n \omega_{X/Y}$ is a reflexive coherent subsheaf of the sheaf
of meromorphic functions of $X$,
and therefore we have
$$
H^0(X, I^n \omega_{X/Y}) = \bigcap_{x \in X: \, \text{ht}(x) = 1}
  \left( I^n \omega_{X/Y} \right)_x.
$$
For $x \in X$ with $\varphi(x) \notin \V(I)$ we have
$$
\eqalignno{
   I \hbox{\Cal O}_{X,x} = & \,  \hbox{\Cal O}_{X,x}, \cr
  \omega_{X/Y,x} = & \,  \hbox{\Cal O}_{X,x}, \cr
}
$$
as $\varphi$ is an isomorphism away from $\V(I)$.
Those $x \in X$ with $\text{ht}(x) = 1$ and
$\varphi(x) \in \V(I)$ correspond to the Rees valuations of $I$,
and thus
  $$
  \eqalignno{
  \adj_X( I^n) = & \, \,  H^0(X, I^n \omega_{X/Y}) \cr
  = & \, \,  \bigcap_{x \in X: \text{ht}(x) = 1}
  \left( I^n  \omega_{X/Y} \right)_x \cr
  = & \, \, \bigcap_{x \in X: \text{ht}(x) = 1, \varphi(x) \in \sV(I)}
  \left( I^n  \omega_{X/Y} \right)_x \, \cap \bigcap_{x \in X: \, 
  \text{ht}(x) = 1, \varphi(x) \notin \sV(I)}
  \hbox{\Cal O}_{X,x} \cr
  \supseteq & \, \, \bigcap_{v \in \sReesval(I)} \{ r \in R: v(r) \geq n \cdot v(I) 
  -   v(J_{R_v/R}) \},
  \cr
}
$$
where we also use,
that $\omega_{R_v/R}$ is an invertible fractional ideal
with inverse $J_{R_v/R}$.
As  $\pi_* \omega_{X/Y} = \hbox{\Cal O}_X$ by~\cite[\S 4]{LipTei}, the converse 
inclusion is obvious, and we conclude that
 $$\adj_X(I^n) = \bigcap_{v \in \sReesval(I)} \, \{ r \in R ~|~ v(r) \ge n \cdot v(I)
    - v(J_{R_v/R})\} $$
Thus (1) is equivalent to 
 $$ \adj_X(I^n) = \adj_Z(I^n) \quad \quad \text{ for all } n \in \bbN $$
for all $f: Z\rightarrow Y$ as above.

First assume (2).
This direction is implicit in~\cite{Lipmanadj}, cf.~\cite{Lipmanadj}, 1.3.2(b). 
Let $f : Z \rightarrow Y$ be as above. By the assumptions we have trivially
  $$ H^0(X,  I^n \pi_* \omega_{Z/Y}) = H^0(X, I^n \omega_{X/Y}) $$
implying by the calculations preceeding ($*$) that $\adj_X(I^n) = \adj_Z(I^n)$ for all positive integers $n$. Thus (1) folllows.

Conversely suppose that (1) holds, i.e., that $\adj(I^n) = \adj_X(I^n)$ for all positive integers $n$. Then by ($*$) we must have that the canonical inclusions
  $$ H^0(X, I^n \pi_* \omega_{Z/Y}) = \adj_Z(I^n)  \hookrightarrow 
    \adj_X(I^n) =  H^0(X, I^n \omega_{X/Y}) $$
are isomorphisms for all positive integers $n$. If $X'$ denotes the blow-up of $I$ on $Y$, then $I \hbox{\Cal O}_{X'}$ is a very ample invertible sheaf on $X'$. As $X/X'$ is finite, $I \hbox{\Cal O}_X$ is an ample invertible sheaf on $X$, and thus the above isomorphisms imply that the canonical inclusion
  $$ \pi_* \omega_{Z/Y} \hookrightarrow \omega_{X/Y} $$
is an isomorphism, i.e., that (2) holds.
\qed

\bigskip

For both conditions in the theorem some examples are known, as we show below.

\smallskip

Recall that two ideals $I, J \subseteq R$ are called projectively equivalent if 
there exist positive integers $i,j$ with $\overline{I^i} = \overline{J^j}$, cf.~\cite{CHRR}.

\cor
\label{\cormonompsrat}
Let $R$ be a regular domain, let $x_1, \ldots,
x_d$ be a generalized regular system of parameters, and let $I$ be
an ideal projectively equivalent to an ideal generated by monomials $\underline 
x^{\underline a_1}, \ldots, \underline x^{\underline a_s}$
in $x_1, \ldots, x_d$.
Then $\adj(I)$ is a monomial ideal in $x_1, \ldots , x_d$, determined by the 
Rees valuations of $I$, and the normalized blow-up of $I$ satisfies condition 
(2) of~\thmadjone.
\endb

\proof 
It remains to note that
  $$ \text{Proj}(\overline{R[It]}) = \text{Proj}(\overline{R[I^it]}) =
   \text{Proj}(\overline{R[\overline{I^i}t]}).  $$
Then the corollary follows from Theorem~{\thmadjmonomval}.
\qed

\medskip

By the work of Lipman and Teissier we also know (2) in some cases.

\cor
\label{\coradjtwo}
Let $(R,\m)$ be a two-dimensional regular domain.
Then for any non-zero ideal $I$,
$$
\adj(I) =
\bigcap_{v \in \sReesval(I)} \{ r \in R ~|~ v(r) \ge v(I) - v(J_{R_v/R})\}.
$$
\endb

\proof
It remains to note that in the two-dimensional case
the normalized blow-up of $I$ has pseudo-rational singularities only
by~\cite{LipTei}, p.~103 and~\cite{LipmanRS}, 1.4.
Thus condition (2) of ~\thmadjone $\, $ is satisfied.
\qed

\medskip

\remark
\label{\rmksubaddtwo} 
In the case of two-dimensional regular rings an elementary direct proof of~\coradjtwo $\, \,$ can be given as well: We may assume that $(R, \m)$ is local with infinite residue field and that $I$ is $\m$-primary. Then it follows 
from~\cite{Hu87} and ~\cite{HunSwa} (see also~\cite{Lipmanadj}) that for a generic 
$x \in \m \setminus \m^2$ the ideals $I$ and $\adj(I)$ are contracted from 
$S := R[\frac {\m}{x}]$ and that $\adj(I) S = \frac {1}{x} \adj(IS)$. From 
this~\coradjtwo $\, $ follows by an easy induction on the multiplicity 
$\text{mult}(I)$ of $I$.

With this line of argument we can also give an easy proof of subadditivity of adjoint ideals in the two-dimensional case. Again we may assume that $(R, \m)$ is local with infinite residue field, and that $I$ and $J$ are $\m$-primary. For a generic $x \in \m \setminus \m^2$ we will have that $I$, $J$, $IJ$, $\adj(I)$, $\adj(J)$, 
$\adj(IJ)$ and $\adj(I)\adj(J)$ are contracted from $S = R[ \frac {\m }{x} ]$. Denoting by $I'$, resp. $J'$, the strict transforms of $I$, resp. $J$, we conclude by the above and by 
induction on $\text{mult}(I) + \text{mult}(J)$:
$$ \eqalignno{
   \adj(IJ)
   & = \adj(IJ) S \cap R \cr
   & = \frac {1}{x} \adj(IJS) \cap R  \cr
   & = x^{\ord(I) + \ord(J)-1} \adj(I'J') \cap R \cr
   & \subseteq x^{\ord(I)-1} \adj(I') \cdot x^{\ord(J)-1} \adj(J') \cap R \cr
   & = \adj(I) \adj(J) S \cap R  \cr
   & = \adj(I) \adj(J).
  } $$
Alternatively, the subadditivity result may be deduced from~\cite{Lipproxy} and~\cite{Lipmanadj}. We note that Tagaki and Watanabe~\cite{TakWat04a} 
proved subadditivity of adjoint ideals more generally, for two-dimensional 
log-terminal singularities. The argument given here does not extend to their situation.

\medskip
\noindent {\tenbf Acknowledgements.}
We thank the referee for pointing out a crucial simplification.

\bigskip\bigskip
\leftline{\bf References}
\bigskip

\font\eightrm=cmr8 \def\rm{\fam0\eightrm}
\font\fiverm=cmr5 \def\rm{\fam0\eightrm}
\font\eightit=cmti8 \def\it{\fam\itfam\eightit}
\font\eightbf=cmbx8 \def\bf{\fam\bffam\eightbf}
\font\eighttt=cmtt8 \def\tt{\fam\ttfam\eighttt}
\textfont0=\eightrm \scriptfont0=\fiverm
\rm
\baselineskip=9.9pt
\parindent=3.6em

\setbox1=\hbox{[999]} 

\newif\ifnumberbibl \numberbibltrue

\newdimen\labelwidth \labelwidth=\wd1 \advance\labelwidth by 2.5em
\ifnumberbibl\advance\labelwidth by -2em\fi

\thmno=0
\def\bitem#1{\advance\thmno by 1%
  \ifnumberbibl
    \immediate\write\isauxout{\noexpand\expandafter\noexpand\gdef\noexpand\csname 
#1\noexpand\endcsname{\the\thmno}}
    \item{\hbox to \labelwidth{\hfil\the\thmno.\ \ }}
  \else
    \immediate\write\isauxout{\noexpand\expandafter\noexpand\gdef\noexpand\csname 
#1\noexpand\endcsname{#1}}
    \item{\hbox to \labelwidth{\hfil#1\ \ }}
  \fi}
\parindent=\labelwidth \advance\parindent by 0.5em 


\bitem{CHRR}
C. Ciuperca, W. Heinzer, L. Ratliff and D. Rush,
Projectively equivalent ideals and Rees valuations,
{\it J. Algebra} {\bf 282} (2004), 140--156.

\bitem{DEL}
J.-P. Demailly, L. Ein and R. Lazarsfeld,
A subadditivity property of multiplier ideals,
{\it Mich. Math. J.} {\bf 48} (2000), 137--156.

\bitem{Howald}
J. A. Howald,
Multiplier ideals of monomial ideals,
{\it Trans. Amer. Math. Soc.} {\bf 353} (2001), 2665--2671.

\bitem{Hu87}
C. Huneke,
Complete ideals in two-dimensional regular local rings, 
{\it Commutative algebra (Berkeley, CA, 1987)} (1989), Math. Sci. Res. 
Inst. Publ., {\bf 15},
Springer, New York, {325--338}.

\bitem{HunSwa}
C. Huneke and I. Swanson,
Cores of ideals in two dimensional regular local rings,
{\it Michigan Math. J.} {\bf 42} (1995), 193--208.

\bitem{Kap}
I. Kaplansky,
$R$-sequences and homological dimension,
{\it Nagoya Math. J.} {\bf 20} (1962), 195--199.

\bitem{KiSt}
K. Kiyek and J. St\"uckrad,
Integral closure of monomial ideals on regular sequences, {\it Rev. Math. 
Iberoamericana} {\bf 19} (2003), 483--508.

\bitem{LipmanRS}
J. Lipman,
Desingularization of two-dimensional schemes, {\it Ann. of Math.}
{\bf 107} (1978), 151--207.

\bitem{Lipmanadj}
J. Lipman,
Adjoints of ideals in regular local rings,
with an appendix by S. D. Cutkosky,
{\it Math. Research Letters} {\bf 1} (1994), 739--755.

\bitem{Lipproxy}
J. Lipman,
Proximity inequalities for complete ideals in two-dimensional regular local rings,
{\it Contemp. Math.} {\bf 159} (1994), 293--306.

\bitem{LipSat}
J. Lipman and A. Sathaye,
Jacobian ideals and a theorem of Brian\c con-Skoda,
{\it Michigan Math. J.} {\bf 28} (1981), {199--222}.

\bitem{LipTei}
J. Lipman and B. Teissier,
Pseudo-local rational rings
and a theorem of Brian\c con-Skoda about integral closures of ideals,
{\it Michigan Math. J.} {\bf 28} (1981), {97--112}.

\bitem{MS63}
H. Muhly and M. Sakuma,
Asymptotic factorization of ideals,
{\it J. London Math. Soc.} {\bf 38} (1963), 341--350.

\bitem{Ra74}
L. J. Ratliff, Jr.,
Locally quasi-unmixed Noetherian rings and ideals of the principal class,
{\it Pacific J. Math.} {\bf 52} (1974), {185--205}.

\bitem{Rees}
D. Rees,
Valuations associated with ideals (II),
{\it J. London Math. Soc.} {\bf 36} (1956), 221--228.

\bitem{SwHu}
I. Swanson and C. Huneke,
{\it Integral Closure of Ideals, Rings, and Modules},
book, to be published by Cambridge University Press, 2006.

\bitem{TakWat04a}
S. Takagi and K.-I. Watanabe,
When does subadditivity theorem for multiplier ideals hold?
{\it Trans. Amer. Math. Soc.} {\bf 356} (2004), 3951--3961.

\bigskip

\noindent
NWF I - Mathematik, Universit\"at Regensburg, 93040 Regensburg,
Germany,

\noindent {\tt Reinhold.Huebl@Mathematik.Uni-Regensburg.de}.

\noindent
Reed College,
3203 SE Woodstock Blvd,
Portland, OR 97202, USA,
\noindent {\tt iswanson@reed.edu}.

\end